\documentclass[leqno,12pt]{article}
\usepackage{amsmath, amssymb}
\usepackage[dvips]{graphicx}
\usepackage{amsthm}

\usepackage{color}%\UTF{0092}\UTF{00C7}\UTF{0089}\UTF{00C1}%
\usepackage{ulem}%\UTF{0092}\UTF{00C7}\UTF{0089}\UTF{00C1}%
\theoremstyle{plain} % イタリック体
\newtheorem{theorem}{\indent\sc Theorem}[section] % 見出しはスモールキャップ
\newtheorem{lemma}[theorem]{\indent\sc Lemma}

\theoremstyle{definition} % ローマン体に変更

\newtheorem{remark}[theorem]{\indent\sc Remark}

\title{\uppercase{Parallel mean curvature tori \\
in $\mathbb{C}P^{2}$ and  $\mathbb{C}H^{2}$} }
\author{ \textsc{Katsuei  Kenmotsu$^{*}$} }

\makeatletter

\@addtoreset{equation}{section}
\makeatother

\date{}

\begin{document}
\maketitle
\footnote{
$^{*}$Partly supported by the Grant-in-Aid for Scientific Research (C), Japan Society for the Promotion of Science.}
\begin{abstract}
 We explicitly determine  tori  that have a parallel mean curvature vector, both  in the complex projective plane  and the complex hyperbolic plane.
\end{abstract}
\section{Introduction}
Parallel mean curvature tori in  real four-dimensional  space forms are  tori with constant mean curvature(cmc) that lie on a totally geodesic hypersurface  or  an umbilical hypersurface of the ambient space, as proved    independently by Hoffman  \cite{hoffman}, Chen \cite{chen1}, Yau \cite{yau}.  The study of such cmc tori in  real three-dimensional space forms is a  topic in the  theory of submanifolds (see, for example,  Wente \cite{wente}, Abresch \cite{abresch}, Pinkall and Sterling \cite{pinksterl}, Bobenko \cite{bobenko}, and
Andrews and Li \cite{andli}).  In the above-mentioned studies, parallel mean curvature tori in the complex two-plane $\mathbb{C}^{2}$ are classified by identifying  $\mathbb{C}^{2}$ with  real Euclidean four-space.  In this paper, we study such tori in  non-flat  two-dimensional complex space forms. The main result of this paper tells us that the non-flat case is much simpler than the flat one. In fact,   Theorem 4.4 of this paper states that a non-zero parallel mean curvature torus in  a non-flat  two-dimensional complex space form is flat and  totally real. Then,  in Theorem 4.5, applying  a result of Ohnita \cite{ohnita} or Urbano \cite{urbano},  we explicitly determine   parallel mean curvature tori  in the complex projective plane $\mathbb{C}P^{2}$  and the complex hyperbolic plane $\mathbb{C}H^{2}$, endowed with their canonical structures of Kaehler surfaces.

 We remark that parallel mean curvature spheres in $\mathbb{C}P^{2}$ and $\mathbb{C}H^{2}$ have been classified  by Hirakawa \cite{hirakawa} and Fetcu \cite{fetcu}. Additional studies regarding  parallel mean curvature surfaces in  various ambient spaces that are neither real nor complex space forms can be found in Torralbo and Urbano \cite{torurb}, Fetcu and Rosenberg \cite{fetros1, fetros2, fetros3}, and Ferreira and Tribuzy \cite{fertri}.

In this paper,  we first study the local structure of  a parallel mean curvature surface  in a  two-dimensional  complex space form using different methods from  those in our previous works  %Kenmotsu and Zhou 
\cite{kenzhou, ken1}.  The new ingredient in this paper is Lemma 2.4, which proves the existence of   special isothermal coordinates fitting to the geometry of a parallel mean curvature vector.  Then, we  determine explicitly the second fundamental form of the surface for these coordinates. By coupling these results with  generalised Hopf differentials found  in Ogata \cite{ogata} and   Fetcu \cite{fetcu}, we prove our main results,  Theorems 4.4 and 4.5.
 % \footnote[0]{Partly supported by  JSPS Grant-in-Aid for Scientific  Research (C-25400062)} 

\section{Preliminaries}
Let $\overline{M}[4\rho]$ be a complex  two-dimensional complex space form with  constant holomorphic sectional curvature $4\rho$. Furthermore, let
$M$ be an oriented and connected real two-dimensional Riemannian manifold with  Gaussian curvature
$K$ and $x: M\longrightarrow
\overline{M}[4\rho]$ be an isometric immersion, with   Kaehler angle  $\alpha$ such that the  mean
curvature vector  $H$  is nonzero and parallel for the normal connection on the normal bundle
of the immersion. We remark that
the case with $H=0$ has been studied in Chern and Wolfson \cite{cherwolf} and  Eschenburg, Guadalupe, and Tribuzy \cite{escgt}.  

Let $M_{0} = \{p \in M\ |\  x \ \mbox{is neither holomorphic nor anti-holomorphic at} \  p \}$.  $M_{0}$ is an open dense subset of $M$. Because all of the  calculations and formulas on $M_{0}$ presented  in Ogata \cite{ogata} are valid until page 400, according to the remark in Hirakawa \cite{hirakawa},  there exists a local
field of unitary coframes
$\{w_{1},w_{2}\}$ on $M_{0}$ such that, by restricting it to $x$, the Riemannian metric $ds^2$
 on $M_{0}$ can be written as  $ds^{2}=\phi\bar{\phi}$, where
$\phi=\cos\alpha/2\cdot\omega_{1}+\sin\alpha/2\cdot\bar{\omega}_{2}$ % (see \cite{ken1})
.  Let $a$ and $c$ be the complex valued functions on $M_{0}$ that determine the second fundamental form of $x$.
 Then,  the Kaehler angle $\alpha$ and the complex 1-form $\phi$  satisfy
\begin{eqnarray}
d\alpha &=& ( a +  b )\phi + (\bar{a} + b )\bar{\phi} , \\
d\phi &=&  (\bar{a} - b)\cot \alpha \cdot \phi \wedge \bar{\phi},
\end{eqnarray}
where $2b = |H| >0$. By (2.4), (2.5), and (2.6) of Ogata \cite{ogata}, we have  
 \begin{eqnarray}
&& K= -4(|a|^{2} - b^{2})+6\rho\cos^{2}\alpha, \\
&& da\wedge\phi= -\left(2a (\bar{a} - b) \cot\alpha  + \frac{3}{2}\rho\sin \alpha \cos \alpha \right)
\phi\wedge\bar{\phi}, \\
&& dc\wedge\bar{\phi}=2c(a - b)\cot\alpha \cdot \phi\wedge\bar{\phi}, \\
&& |c|^{2} = |a|^{2} + \frac{\rho}{2}(-2 + 3\sin^{2}\alpha),
\end{eqnarray}
where (2.3) is the Gauss equation, (2.4) and (2.5) are the Codazzi-Mainardi equations, and (2.6) is the Ricci equation of $x$.
\begin{remark}  (1)\ The immersion $x$ is holomorphic (resp. anti-holomorphic) at $p \in M$ if and only if $\alpha=0\  (\mbox{resp.}\  \alpha = \pi)$ at $p$. Hence, it holds that $\sin \alpha \neq 0$  on  $M_{0}$.

(2)\  The unitary coframes $\{w_{1},w_{2}\}$ used in $(2.1) \sim (2.5)$ are uniquely determined up to the orientations of both $\overline{M}[4\rho]$ and $M_{0}$ (see \cite{ogata}), hence the complex one-form $\phi$ on $M_{0}$ is unique up to its sign and  conjugacy.
\end{remark}

Let us take an isothermal coordinate for the Riemannian metric $\phi \bar{\phi}$ on $M_{0}$ that makes $M_{0}$ a Riemann surface with a local complex coordinate $z$, and  put
$$
\Phi _{1} = (8ba - 3\rho\sin^{2}\alpha)\phi^2, \quad \Phi_{2} =\bar{c}\phi^2.
$$
We will use the following facts that were proved in  Ogata \cite{ogata} and Fetcu \cite{fetcu}.
\begin{lemma} The quadratic forms $\Phi_{1}$ and $\Phi_{2}$ on $M_{0}$ are holomorphic.
\end{lemma}
For completeness, we will present a proof of Lemma 2.2. The quadratic forms $Q$ and $Q'$ on $M$ in Fetcu \cite{fetcu} are written in our terminology as
\begin{equation}
Q=(8b(\bar{c}+a) - 3 \rho \sin^2\alpha)\phi^2, \quad Q'= (8b(\bar{c}-a) + 3 \rho \sin^2\alpha) \phi^2.
\end{equation}
Let $\phi = \lambda dz \ (\lambda \neq 0)$. By using $(2.1) \sim (2.5)$, we see that
$$
d(\lambda^2 (8b(\bar{c}+a) - 3 \rho \sin^2\alpha))\wedge \phi =0, \   d(\lambda^2 (8b(\bar{c}-a) + 3 \rho \sin^2\alpha))\wedge \phi =0,
$$
 showing that  both forms of (2.7) are holomorphic on $M_{0}$. By taking both the addition and the subtraction of  the two forms of (2.7), we  proved Lemma 2.2.
\vspace{0.5cm}
 
The local structure of the  immersion satisfying $a =\bar{a}$ on $M_{0}$ has been determined  in    Kenmotsu and Zhou \cite{kenzhou} and Kenmotsu \cite{ken1}.  We remark that the classification table of those surfaces can also be found  in   Hirakawa \cite{hirakawa}.  By applying these results, we obtain the following lemma.
\begin{lemma} Let $\rho \neq 0$ and let $x: M  \longrightarrow
\overline{M}[4\rho]$ be an isometric immersion with a  non-zero parallel   mean
curvature vector.  If  $a=\bar{a}$ on $M_{0}$, then either  $K \equiv 0$ on $M$ and $x$ is totally real,  or $K\leq -2b^2$ on $M$ and $x$ is not totally real.
\end{lemma}
Proof. If $x$ is totally real, then $\alpha= \pi/2$, so $K = 0$ on $M_{0}$ by (2.1) and (2.3). By the continuity of $K$, we have also $K = 0$ on $\overline{M_{0}}=M$.
Suppose that $x$ is not totally real. When $\rho \neq 0$,  we have that  $\rho = - 3b^2$ \cite{kenzhou, ken1},  and  by   Theorem 2.1  (iii) in Hirakawa \cite{hirakawa} we have $K \leq  -2b^2$. In fact, we see  $K= \ \mbox{constant} \ = -2b^2$ on $M_{0}$  if $\alpha$ is constant,  and  
$K=-2( 2(a+b )^2 +b^2)  \leq -2b^{2}$ on $M_{0}$ if $\alpha$ is not constant. By the continuity of $K$, we proved Lemma 2.3.
\vspace{0.5cm}

From this point on, we will  study  an immersion satisfying $ d\alpha \neq 0$ and  $a \neq \bar{a}$  at a point of $M_{0}$.
An immersion $x$ is called  a general type if it satisfies these two conditions  on $M_{0}$.

  We proved in \cite{ken1} that if $x$ is of a general type, then
$a$ is a function of $\alpha$, say $a=a(\alpha)$,  satisfying the first order ordinary complex differential equation
\begin{equation}
 \frac{da}{d\alpha} = \frac{\cot \alpha}{\overline{a(\alpha)} +b}\left( - 2ba(\alpha) +  2|a(\alpha)|^{2}
 + \frac{3\rho}{2}\sin^{2}\alpha \right), \ (a+b \neq 0)
 \end{equation}
and  the Kaehler angle  $\alpha$ satisfies the second order partial differential equation 
\begin{equation}
\alpha_{z\bar{z}} - F(\alpha)\alpha_{z}\alpha_{\bar{z}} =0,
\end{equation}
where
\begin{equation}
F(\alpha) = \frac{((a(\alpha) - b)(\overline{a(\alpha)} - b)
+ 3\rho/2 \sin^{2} \alpha)}{(a(\alpha) + b)(\overline{a(\alpha)} + b)} \cot \alpha .
\end{equation}
The following Lemma plays  a fundamental role in  this paper.
\begin{lemma} Suppose that $x$ is of  a general type.  Then, there is a  coordinate transformation $w=w(z)$ on the Riemann surface $M_{0}$ such that  in $\phi = \mu dw$, $\mu$ is a complex valued function of  the single real variable  $(w+\bar{w})/2$.
\end{lemma}
Proof. By (2.9), we have that $\alpha_{z\bar{z}}/\alpha_{z} - F(\alpha)\alpha_{\bar{z}} =0$, which
 implies that the function  $\log \alpha_{z} - \int F(\alpha)d\alpha$ is holomorphic for $z$. Hence, there exists a holomorphic function $G(z)$ such that
$\alpha_{z} = G(z) \exp(\int F(\alpha)d\alpha)$. We set $w= \int G(z)dz$. Then, it follows from (2.1)  that
$
\phi = \lambda dz = \alpha_{z}/(a(\alpha)+b)dz=\mu dw,
$
where we set
\begin{equation}
\mu =  \frac{\exp \left(\int F(\alpha)d\alpha \right)}{a(\alpha)+b}.
\end{equation}
We remark that $\mu$ is a complex valued function of $\alpha$. Now, we prove that $\alpha$ is a function of $u$, say $\alpha = \alpha(u)$, where 
we set $w=u +i v, (u, v \in \mathbb{R}).$ In fact, by (2.1) and (2.11) we can see that
\begin{equation}
d\alpha =(a+b)\mu dw + (\bar{a}+b)\bar{\mu}d\bar{w}
= 2 \exp \left(\int F(\alpha)d\alpha \right) du,
\end{equation}
because $F(\alpha)$ is real valued,  which proves Lemma 2.4.
\vspace{0.3cm}

We note that the Kaehler angle  $\alpha$ and the complex valued function $\mu$ are now   functions of a single real variable  $u$ for  coordinates $(u,v)$. It follows from (2.8), (2.10), and (2.11) that
\begin{equation}
\frac{d\log \mu}{d\alpha} = - \frac{(\overline{a(\alpha)}-b)}{\overline{a(\alpha)}+b}\cot \alpha.
\end{equation}
We can determine $c$ as follows.
\begin{lemma} Suppose that $x$ is of a general type. Then, there is a    real number $k_{1}$ such   that
\begin{equation}
c= \left(|a(\alpha)|^{2} + \frac{\rho}{2}\left(-2 + 3 \sin^{2}\alpha \right) \right)^{1/2} \frac{(\overline{a(\alpha)}+b)}{a(\alpha)+b} e^{-i k_{1}v}.
\end{equation}
\end{lemma}
 Proof. Set $a(\alpha) +b= |a(\alpha)+b|\exp i \theta(\alpha)$ and $c =|c|\exp i\nu(u,v)$, where 
 $\theta(\alpha)$ and $\nu(u,v)$ are real valued functions of $\alpha$ and $u$ and $v$ respectively.  By (2.11), the absolute value and the argument of $\mu$ are  functions of $\alpha$ and $-\theta(\alpha) \ ( \mbox{mod} \  2\pi)$, respectively.  Hence,
we have that 
$$
\mu^{2}\bar{c} = |\mu|^{2}(\alpha) |c|(\alpha) \exp i(-2\theta(\alpha) - \nu(u,v)).
$$
We note that $\theta = \theta(\alpha(u))$ is a function of $u$ only. Hence, the absolute value of $\mu^{2}\bar{c}$ 
is a function of $u$ only. Together with Lemma 2.2, this  implies that
\begin{equation*}
 \frac{d}{d\alpha}\left(|\mu|^2|c| \right)\frac{d\alpha}{du} + |\mu|^2|c|\frac{\partial \nu}{\partial v} =0, \quad
 2 \theta'(\alpha) \frac{d\alpha}{du} + \frac{\partial \nu}{\partial u} =0.
\end{equation*}
By the second equation above, we see that $ \nu(u,v) = - 2\theta(\alpha(u)) + \xi(v)$ for some real valued function $\xi$ of $v$.
 By  inserting this to the first equation, we find that there is a real number $k_{1}$ satisfying 
\begin{equation}
\frac{d}{du}\left(\log \left (|\mu|^2|c|(\alpha(u)\right)\right) = k_{1},  \quad
\xi'(v) = - k_{1}.  
\end{equation}
By the second formula of (2.15), we have that $\nu(u,v)= -2\theta(\alpha(u))  - k_{1}v + k_{2}$ for some $k_{2} \in \mathbb{R}$. Because $c$ is defined up to a complex factor of norm $1$, we may assume that $k_{2}=0$. Then, by the fact that $(\bar{a}+b)/(a+b) = \exp(-2i\theta)$,  Lemma  2.5 is proved.
\vspace{0.5cm}

The constant $k_{1}$ admits the following expression. 
\begin{lemma} Suppose that $x$ is of a general type. Then, we have
\begin{equation}
2k_{1}
= \rho\mu\frac{\left( 8 |a|^2 + 9b(a+\bar{a}) \sin^2 \alpha - 8 b^2 + 18 b^2\sin^2\alpha\right)}
{ (\bar{a}+b) \left(|a|^{2} + \rho/2(-2 + 3\sin^{2}\alpha)\right)}  \cot \alpha .
\end{equation}
\end{lemma}
Proof.   By the first equation of  (2.15) and (2.6), we have
$$
%k_{1} &=& \frac{d]{d\alpha}(\log \mu + \log \bar{\mu} + \frac{1}{2} \log {c|^{2}) \frac{d\alpha}{du} \\
%k_{1}&=& \frac{d\alpha}{d u}\left( \frac{d \log \mu}{d\alpha} +  \frac{d \log \bar{\mu}}{d\alpha} + \frac{1}{2} \frac{1}{|c|^{2}} \frac{d|c|^{2}}{d\alpha} \right) \ \mbox{by}\ (2.15) \\
k_{1}= \frac{d\alpha}{d u}\left( \frac{d \log \mu}{d\alpha} +  \frac{d \log \bar{\mu}}{d\alpha} + \frac{1}{2} \frac{1}{|c|^{2}}\left( \frac{d a}{d\alpha} \overline{a(\alpha)} + a(\alpha) \frac{d \bar{a}}{d\alpha}  + 3 \rho \sin \alpha \cos \alpha\right )   \right).
  %&=& 2 e^{\int F(\alpha) d\alpha} \left(-\frac{(\bar{a}-b)}{\bar{a}+b} \cot \alpha - \frac{(a-b)}{a+b} \cot \alpha +    \frac{1}{2} \frac{1}{|c|^{2}}\left( \frac{d|a(\alpha)|^{2}}{d\alpha}  + 3 \rho \sin \alpha \cos \alpha \right) \right).
$$
%&=& 2 e^{\int F(\alpha) d\alpha} (-\frac{(\bar{a}-b)}{\bar{a}+b} \cot \alpha - \frac{(a-b)}{a+b} \cot \alpha   +  \frac{1}{2} \frac{1}{|c|^{2}}( \frac{d a}{d\alpha} \bar{a} + a \frac{d \bar{a}}{d\alpha}  + 3 \rho \sin \alpha \cos \alpha ) )  
%\end{eqnarray*}
It follows from (2.8), (2.11), (2.12), and (2.13) that  (2.16) holds, proving Lemma 2.6.
\vspace{0.2cm}

Lemma 2.2 can be sharpened for the case of a general type.
\begin{lemma} Suppose that $x$ is of a general type. Then, there exist some constants $c_{1}, c_{2} \in \mathbb{C}$ such that
\begin{equation}
\mu^{2}(8ba - 3 \rho \sin^{2}\alpha) =c_{1}, \quad \mu^{2}\bar{c} = c_{2}e^{k_{1}\omega}.
\end{equation}
\end{lemma}
Proof.   By Lemma 2.2, $\mu^{2}(8ba - 3 \rho \sin^{2}\alpha)$ is a holomorphic function for $\omega$. Since we know $\mu = \mu(\alpha), \ a = a(\alpha)$ and $\alpha = \alpha(u)$, this is a function of $u$ only, which implies the first formula of (2.17).  By (2.11) and (2.14), the argument of $\mu^{2}\bar{c}$ is $k_{1}v$ and by the first formula of (2.15), its absolute value is $c_{2}e^{k_{1}u}$,  which proves Lemma 2.7.
\section{The case with $k_{1}=0$}
In this section, we study the case with  $k_{1}=0$.
The main result of this section  is used to prove the main theorems of this paper.  Let us recall  assumptions to be applied in this section. These are that $x$ is of a general type,  
$
\rho  \neq  0, \  a+b  \neq  0, \  c  \neq  0, \  \sin \alpha  \neq  0\ \mbox{(by Remark 2.1)},\ \mbox{and} \   k_{1}=0.
$

By (2.16),  we have 
\begin{equation}
 8 |a(\alpha)|^2 + 9b(a(\alpha)+\overline{a(\alpha)}) \sin^2 \alpha - 8 b^2 + 18 b^2\sin^2\alpha=0.
\end{equation}

Moreover, we prove
\begin{lemma} If $x$ is of a general type with $k_{1}=0$, then it holds that
$$
\cot \alpha \left(\sin^{2}\alpha - \frac{8}{9} \right) \left(|a|^2 - b^2 \right) \left(\rho + 3 b^2 \right) =0.
$$
\end{lemma}
Proof.  Let us take the derivative of (3.1) for $\alpha$. Using (3.1) and $d\bar{a}/d\alpha = \overline{da/d\alpha}$, we have
\begin{equation}
 8\frac{d}{d\alpha}|a(\alpha)|^{2} + 9 b \left(\frac{da}{d\alpha}+\overline{\frac{da}{d\alpha}} \right)\sin^{2}\alpha - 16 \cot \alpha \left(|a|^{2}-b^2 \right)=0.
\end{equation}
By (2.8) and (3.1), we have
\begin{eqnarray}
&& \frac{d}{d\alpha}|a(\alpha)|^{2} = \frac{\cot \alpha}{|a+b|^{2}} \left(|a|^{2}-b^{2} \right) \left(4 |a|^{2} - \frac{4}{3} \rho + 3 \rho \sin^{2}\alpha\right ),  \\
&& \frac{da}{d\alpha}+\overline{\frac{da}{d\alpha}}  = 
 4b \frac{\cot \alpha}{|a+b|^{2}} \left(|a|^{2}-b^2 \right)   \nonumber \\
&&\qquad \qquad  + \frac{\cot \alpha}{|a+b|^{2}} (a+\overline{a}+2b)\left(-2b(a+\bar{a})+2b^2 + 2|a|^{2}+\frac{3}{2}\rho \sin^{2}\alpha \right). \nonumber
\end{eqnarray}
By (3.1), we have $9b(a+\bar{a}+2b) \sin^{2}\alpha = - 8 (|a|^{2}-b^2)$. Coupling this with the second formula above, we have
\begin{eqnarray*}
  &&    9b\left( \frac{da}{d\alpha}+\overline{\frac{da}{d\alpha}}
\right )\sin^{2}\alpha   \\
&&    =  \frac{\cot \alpha}{|a+b|^{2}}\left(|a|^{2}-b^2 \right) \left(-16|a|^{2}+16b(a+\bar{a})-16b^{2}+36b^{2}\sin^{2}\alpha-12 \rho \sin^{2}\alpha \right). \nonumber
\end{eqnarray*}
Lemma 3.1 is proved by this formula above, (3.2), and (3.3).

Since $x$ is of a general type, $\alpha$ is not constant, hence we may assume $\cot \alpha \neq 0$, and $\sin^{2}\alpha - 8/9 \neq 0$. If $|a(\alpha)|^2=b^2$ on $M_{0}$, then we have by (3.1) that $a(\alpha) + \overline{a(\alpha)}=-2b$, and so $(a(\alpha)+b)(\overline{a(\alpha)}+b)=0$ on $M_{0}$,  which contradicts the assumptions applied in this section. Hence, Lemma 3.1 implies that 
\begin{equation}
\rho=-3b^2.
\end{equation}
By considering (2.14) with $k_{1}=0$,  we note that $c$ is a function of $u$ only, because   $\alpha =\alpha(u)$.  In order to determine $a(\alpha)$ explicitly, we 
set 
\begin{equation}
 \tau =\frac{a(\alpha) - b}{a(\alpha) + b}.
\end{equation}
Then, (2.8) and (3.1)  are  transformed to
%\begin{equation}
$$
\frac{d\tau}{d\alpha} -(1+\tau)(1-\tau)\bar{\tau}\cot \alpha +\frac{9}{8}(1-\tau)^{2}(1-\bar{\tau})\sin \alpha \cos \alpha =0,
$$
%\end{equation}
and
$
18\sin^{2}\alpha + (8-9\sin^{2}\alpha)(\tau+\bar{\tau}) =0,
$
respectively. Hence, the real part of $\tau$ is uniquely determined, and the imaginary part of $\tau$, say $y(\alpha)$, satisfies the 
first order ordinary differential equation
\begin{equation}
\frac{dy^{2}}{d\alpha} + 4\cot \alpha\frac{(4-9\sin^{2}\alpha)}{8-9\sin^{2}\alpha}   y^{2} +\cot \alpha \frac{(8-9\sin^{2}\alpha)}{4}  y^{4} =0.
\end{equation}
 The solution of (3.6) is given by
\begin{equation}
y^{2}(\alpha) = \frac{8c_{3}}{(8-9\sin^{2}\alpha)(-c_{3} + \sin^{2}\alpha)} \geq 0, \quad (c_{3} \in \mathbb{R}),
\end{equation}
%w= \frac{-9\sin^{2}\alpha}{8 - 9\sin^{2}\alpha} 
%+ i\frac{2\sqrt{2}k_{3}}{\sqrt{(8 - 9\sin^{2}\alpha)(- k_{3}^{2} + \sin^{2}\alpha)}},  
%\end{equation}
where the integral constant $c_{3}$ satisfies $c_{3} \neq 0$ and $8-9c_{3} \neq 0$. In fact, if $c_{3}=0$, then $y(\alpha)$  identically vanishes, hence $a =\bar{a}$ on $M_{0}$, which contradicts the assumption $a \neq \bar{a}$. If $8-9c_{3} = 0$, then $y^{2} <0$, giving a contradiction. Since the right hand of (3.7) is non-negative, it holds that 
\begin{eqnarray}
\left\{
\begin{array}{rl}
&  0 < c_{3} < \sin^{2}\alpha < \frac{8}{9}, \ \mbox{or} \ \  \frac{8}{9} < \sin^{2}\alpha < c_{3}, \  \mbox{for} \  c_{3} >0,  \\
& \frac{8}{9} < \sin^{2}\alpha \leq  1,  \  \mbox{for} \  c_{3} <0.
\end{array} \right.
\end{eqnarray}
%cf. 2015-02-28.nb, 2015-02-27.nb, 2015-02-27b.nb  Lemmas5(2.15) and 6(2.19),(2.20).nb
Later, we will use the  explicit formula  of $a$.  By (3.5) and (3.7),
the real part of $a(\alpha)$, $\Re a(\alpha)$, is given by
\begin{equation}
\Re a(\alpha) = \frac{b}{8-9c_{3}} \frac{(-16c_{3} + (8+27c_{3})\sin^{2}\alpha - 18 \sin^{4}\alpha)}{\sin^{2}\alpha},
\end{equation}
and also, 
the imaginary part of $a(\alpha)$, $\Im a(\alpha)$,  is given by, when $c_{3}>0$,
\begin{equation}
  \Im a(\alpha) = \frac {b \sqrt{c_{3}}}{\sqrt{2}(8-9c_{3})}\frac{ (8-9\sin^{2}\alpha)\sqrt{(8-9\sin^{2}\alpha)(-c_{3}+\sin^{2}\alpha)}}{\sin^{2}\alpha},  
\end{equation}
and  when $c_{3}<0$, setting $c_{4} =-c_{3} >0$,
\begin{equation}
 \Im a(\alpha) = \frac{b \sqrt{c_{4}}}{\sqrt{2}(8 + 9c_{4})}\frac{ (9\sin^{2}\alpha -8)\sqrt{(9\sin^{2}\alpha-8) (c_{4} +\sin^{2}\alpha)}}{ \sin^{2}\alpha}, 
%\left( \frac{8}{9} < \sin^{2}\alpha <1 \right) .
\end{equation}
where $\sin \alpha$ satisfies (3.8). We note that  by (3.7), (3.10) and (3.11), those intervals in   (3.8) are the maximal ones of $\sin^{2}\alpha$ on which $\Im a(\alpha)$ does not vanish.
We compute the Gaussian curvature of $M_{0}$ as follows: By (3.1) and (3.5), we have
%\begin{eqnarray*}
$$
|a(\alpha)|^{2}=b^2\frac{1+(\tau+\bar{\tau}) +|\tau|^{2}}{1-(\tau+\bar{\tau}) +|\tau|^{2}},\quad
\tau= -\frac{9 \sin^{2}\alpha}{8-9\sin^{2}\alpha} + i y(\alpha)
%\end{eqnarray*}
$$
By these formulas above, (2.3) and  (3.7), we have
\begin{equation}
K=\frac{-2 b^{2}}{(8-9c_{3})}\left( \left(9\sin^{2}\alpha - 8  \right)^{2} + (8-9c_{3}) \right),
%K= \frac{-2b^2}{8-9k_{3}^2}\left( (8 - 9 \sin^{2}\alpha)^2 + 8-9k_{3}^2 \right).
\end{equation}
which proves the following Lemma 3.2
\begin{lemma} If $8-9c_{3} >0$, then the Gaussian curvature $K$ is upper bounded by a negative constant on $M_{0}$, that is,  it holds that $K \leq   -2b^{2} <0$ on $M$.
\end{lemma}
Hence, we have proved part of the following theorem
\begin{theorem}  Let $\rho \neq 0$, and let $x: M\longrightarrow
\overline{M}[4\rho]$ be an isometric immersion with a non-zero  parallel  mean
curvature vector of a general type. If  $k_{1}=0$,  then $\rho=-3b^2$, and the first and second fundamental forms of $x$ are explicitly determined by $(2.11), (2.14), (3.9), (3.10), (3.11)$, and a real number $c_{3}$. In particular, the Gaussian curvature $K$ is bounded above by $-2b^{2}$ when $8-9c_{3}>0$.
Conversely, for   $b>0$,  there is a one-parameter family of parallel  mean curvature  immersions of a general type with $k_{1}=0$  from a simply connected domain $D$ in $\mathbb{R}^{2}$   into $\overline{M}[-12b^2]$,  with $|H| =2b$.
\end{theorem}
Proof of the converse:  For given  $b>0$ and  $c_{3} \in \mathbb{R}$, we can find a simply connected domain $D$ with coordinates $(u,v)$, a real valued function $\alpha =\alpha(u)$, complex valued functions $a=a(\alpha)$ and $c=c(\alpha)$, and a complex one form $\phi$ on $D$, by considering (3.7), (3.9), and (2.14) with $k_{1}=0$. Then, we can  prove that these functions satisfy the structure equations $(2.1) \sim (2.6)$, and that $c_{3}$ is the parameter of the  family. We omit a detailed computation of the proof, because the converse is not used  in this paper.
\begin{remark}
The case with $k_{1} \neq 0$ remains to be  studied. This will be addressed in future work, because the study of this case is not necessary for the present paper.
\end{remark}
\section{ Parallel mean curvature tori }
In this section, we suppose that $M$ is homeomorphic to a torus, and we  denote it  by $\mathbb{T}^{2}$.  Let $\mathbb{T}_{0}^{2} = \{ p \in \mathbb{T}^{2}\ |\  x \ \mbox{is neither holomorphic nor anti-holomorphic at}\  p \}$. First, we study the case of $a=\bar{a}$ on $\mathbb{T}_{0}^{2}$.
\begin{lemma} Let $\rho \neq 0$ and let $x: \mathbb{T}^{2} \longrightarrow \overline{M}[4\rho]$ be an isometric immersion  with a  non-zero parallel mean curvature vector. Suppose that $x$ satisfies $a=\bar{a}$ on $\mathbb{T}_{0}^{2}$. Then,  $K = 0$ on $\mathbb{T}^{2}$ and $x$ is totally real.
\end{lemma}
Proof. By Lemma 2.3, if $x$ is not totally real, then $K \leq -2b^{2}$ on $\mathbb{T}^{2}$. However this presents a contradiction by the Gauss-Bonnet Theorem, proving Lemma 4.1.

From now on, we suppose that $a \neq \bar{a}$ at some point of $\mathbb{T}_{0}^{2}$. In particular, we have $8ba - 3 \rho \sin^{2}\alpha \neq 0$ at the point of $\mathbb{T}_{0}^{2}$, hence by Lemma 2.7  $c_{1} \neq 0$ and $\Phi_{1}$ is a non-zero quadratic form on $\mathbb{T}_{0}^{2}$, hence $Q-Q'$ is a non-zero holomorphic quadratic form on $\mathbb{T}^{2}$ by  2.7 and Fetcu \cite{fetcu}. 
By  Riemann-Roch's Theorem, the dimension of the vector space over $\mathbb{C}$ of holomorphic quadratic forms on a torus is one, and $Q-Q'$ is a base of the vector space. Hence, the holomporphic quadratic form $Q+Q'$ on $\mathbb{T}^{2}$ is a  constant multiple of $Q-Q'$. By Lemma 2.7, we have
\begin{equation}
\mu^2 \bar{c} = \mbox{constant},   \ \mbox{and}  \ k_{1}=0.
\end{equation}
We remark that $c$ in the above formula does not vanish on $\mathbb{T}_{0}^{2}$. In fact, if $c=0$ at a point of  $\mathbb{T}_{0}^{2}$, then it must identically vanish  by (4.1).  Such surfaces  are classified by Hirakawa \cite{hirakawa}. We apply Proposition 3.5 of \cite{hirakawa} to our situation that is $a \neq \bar{a}$ at a point of  
$\mathbb{T}_{0}^{2}$, and get  $K=-2b^{2}$ on  $\mathbb{T}_{0}^{2}$. Hence we have that $K= \mbox{constant} = -2b^{2}$ on  $\mathbb{T}^{2}$,  which presents a contradiction by the Gauss-Bonnet Theorem because of $b \neq 0$. By taking the ratio of the two formulas of (2.17), we have that  for some $\gamma (\neq 0) \in \mathbb{C}$,
\begin{equation}
8ba  - 3\rho \sin^2 \alpha  - b\gamma \bar{c}=0 \quad \mbox{on}\ \mathbb{T}_{0}^{2}.
\end{equation}
We have
\begin{lemma} If $a\neq  \bar{a}$ at a point of $\mathbb{T}_{0}^{2}$, then $|\gamma|^{2}   = 2(8-9c_{3}) >0$.
\end{lemma}
Proof.  Since we have $k_{1}=0$ by (4.1),  we can apply  the results of  Section 3 for $\mathbb{T}_{0}^{2}$.  Then, by (3.4) and  (4.2) we have that
%\begin{equation}
$8 a(\alpha) + 9b \sin^2 \alpha  = \gamma \bar{c}$  on $\mathbb{T}_{0}^{2}$.
%\end{equation}
We take the absolute value of the formula above and then we use (2.6) and (3.1) to get
$$
 4b (8-9 \sin^{2}\alpha)^{2} = |\gamma|^{2}(-9 \sin^{2}\alpha \Re a(\alpha) + 16 b
 -27b \sin^{2}\alpha),
$$
 and by  (3.9)  we  get
$(8-9\sin^{2}\alpha)^{2}( |\gamma|^{2} - 2(8-9c_{3})) =0$. If $|\gamma|^{2} - 2(8-9c_{3})\neq 0$, then $8-9\sin^{2}\alpha=0$ on $\mathbb{T}^{2}$, which implies 
$\alpha=\mbox{constant}$ on $\mathbb{T}^{2}$, so $a=\bar{a}$ on $\mathbb{T}_{0}^{2}$, giving a contradiction. We proved Lemma 4.2.
\begin{lemma} It holds that $a=\bar{a}$ on $\mathbb{T}_{0}^{2}$.
\end{lemma} 
Proof. Let $U = \{ p \in \mathbb{T}_{0}^{2}\ | \  a \neq \bar{a}\  \mbox{at} \ p \}$. Suppose that $U$ is non-empty.  
We note that, by Theorem 3.3 and Lemma 4.2, $K$ is upper bounded by a negative constant on $U$. If $U=\mathbb{T}_{0}^{2}$, then we have a contradiction by the Gauss-Bonnet Theorem. From now on, we consider the case that there is a point $p_{0} \in \mathbb{T}_{0}^{2}$ with $a=\bar{a}$ at $p_{0}$. We remark $c_{3} \neq 0$ because of (3.7).
For $c_{3}>0$, by (3.10),  $\Im a$ tends to zero if and only if $\sin^{2}\alpha$ tends to $c_{3}$ or $8/9$. Therefore, for the case of $c_{3} >0$,  $(c_{3},8/9)$ is the maximal interval of $\sin^{2}\alpha$ satisfying $a\neq \bar{a}$. For the case of $c_{3}<0$, the similar consideration proves that  $(8/9,1)$ is the maximal interval of $\sin^{2}\alpha$ satisfying  $a\neq \bar{a}$.
% Furthermore, we have  by (3.5) that $K \leq -2b^2$ on $U$, proving Lemma 4.2.
We have $ 0 <\sin^{2}\alpha(p_{0}) \leq c_{3}$ or $8/9 \leq \sin^{2}\alpha(p_{0})$. If $\sin^{2}\alpha(p_{0}) < c_{3}$, then there is an open neighborhood $W$ of $p_{0}$ such that $\sin^{2}\alpha <c_{3}$ on $W$. Since we have $a=\bar{a}$ on $W$, it holds that $K \leq -2b^{2}$ at $p_{0}$ by Lemma 2.3. For the  case of $8/9 < \sin^{2}\alpha(p_{0})$,  the similar consideration as before implies also $K \leq -2b^{2}$ on the point of $a=\bar{a}$. Therefore,  for those points $p \in \mathbb{T}_{0}^{2}$ satisfying $a(p)=\bar{a}(p)$, $K$ is also upper bounded by a negative constant. We proved that $K$ is upper bounded above by a negative constant on both sets $U$ and $\mathbb{T}_{0}^{2} \setminus U$,  which contradicts the Gauss-Bonnet Theorem by taking the integral of $K$ over $\mathbb{T}^{2}$. Hence,  we proved Lemma 4.3. 

By Lemmas 4.1, 4.2 and 4.3, we  proved the following main result of this paper:
\begin{theorem}  Let $\rho \neq 0$, and
let  $x: \mathbb{T}^{2} \longrightarrow
\overline{M}[4\rho]$ be an isometric immersion  with a  non-zero parallel mean curvature vector. Then,  the Gaussian curvature of  $\mathbb{T}^{2}$ vanishes identically and $x$ is totally real.
\end{theorem}
Let $\mathbb{C}P^{n}$ and $\mathbb{C}H^{n}$ be the complex projective space  and  complex hyperbolic space, respectively,  endowed with  Kaehler metrics  of constant holomorphic sectional curvature.   In Ohnita \cite{ohnita}, and also independently in Urbano \cite{urbano}, $n$-dimensional totally real submanifolds with non-negative sectional curvature in  $\mathbb{C}P^{n}$ and  $\mathbb{C}H^{n}$ have been classified in the context of the theory of symmetric spaces. For the flat case, these immersions are explicitly described in  Dajczer and Tojeiro  $\cite{dato}$ when the ambient space is $\mathbb{C}P^{n}$, and  Hirakawa $\cite{hirakawa2}$ when the ambient space is $\mathbb{C}H^{2}$.
By combining  Theorem 4.4 with those results, we  determined  tori with non-zero parallel mean curvature in two-dimensional non-flat complex space forms. For completeness, we state the following as a theorem
\begin{theorem}  Let $\mathbb{T}^{2}$ be a real two-dimensional compact orientable  Riemannian manifold with  genus one.  Then, an isometric immersion from $\mathbb{T}^{2}$ into  $\mathbb{C}P^{2}$ or $\mathbb{C}H^{2}$ has a  non-zero parallel mean curvature vector  if and only if the image  is a totally real flat torus in $\mathbb{C}P^{2}$ or $\mathbb{C}H^{2}$, respectively. 
\end{theorem}

\medskip
\begin{flushleft}

 Katsuei  Kenmotsu \\
Mathematical Institute,  Tohoku University  \\
980-8578 \quad  Sendai, Japan \\
email:  kenmotsu@math.tohoku.ac.jp
\end{flushleft}
\end{document}